\documentclass[11pt]{article}
\usepackage{amssymb,epsfig}
\usepackage[all,2cell,ps]{xy}

\newtheorem{theorem}{Theorem}[section]
\newtheorem{proposition}[theorem]{Proposition}
\newtheorem{lemma}[theorem]{Lemma}
\newtheorem{corollary}[theorem]{Corollary}

\newtheorem{problem}[theorem]{Problem}
\newtheorem{example}[theorem]{Example}

\newtheorem{construction}[theorem]{Construction}

\newenvironment{proof}{{\noindent \sc Proof.}}{\hfill $\Qed$}

\newcommand{\Qed}{\rule{2.5mm}{3mm}}

\newcommand{\Aut}{\hbox{{\rm Aut}}\,}
\newcommand{\Cay}{\hbox{{\rm Cay}}\,}

\newcommand{\fib}{\hbox{{\rm fib}}}

\newcommand{\CDC}{\hbox{{\rm CDC}}}

\newcommand{\F}{\mathrm{F}}

\newcommand{\ZZ}{\mathbb{Z}}

\newcommand{\CT}{\hbox{{\rm CT}}}

\newcommand{\tX}{\tilde{X}}
\newcommand{\tG}{\tilde{G}}
\newcommand{\tH}{\tilde{H}}
\newcommand{\bG}{\bar{G}}
\newcommand{\bH}{\bar{H}}

\newcommand{\bg}{\bar{g}}

\newcommand{\p}{\wp}

\newcommand{\C}{{\cal{C}}}

\newcommand{\D}{D}

\newcommand{\la}{\langle}
\newcommand{\ra}{\rangle}

\def\a{\alpha} \def\b{\beta}

 \def\D{\Delta}   
 \def\O{\Omega}

 \def\tX{\widetilde X}

\def\mz{{\mathbb Z}}

\def\Aut{\hbox{\rm Aut}}

\def\mod{\hbox{\rm mod}\;}

\def\Cay{\hbox{\rm Cay}}

\def\Cos{\hbox{\rm Cos}}

\newcommand{\bW}{\bar{W}}

\begin{document}


\begin{center}
{\bf\large ON 2-FOLD COVERS OF GRAPHS}
\end{center}

\bigskip

 Yan-Quan Feng{\small$^{a,}$}\footnotemark, \,
 Klavdija Kutnar{\small$^{b,}$}\footnotemark, \, 
 Aleksander Malni\v c{\small$^{c,}$}\addtocounter{footnote}{-1}\footnotemark$^,$*  \, and
 Dragan Maru\v si\v c{\small$^{b,c,}$}\addtocounter{footnote}{-1}\footnotemark\\

\begin{center}
{\it {\small  $^a$Beijing Jiaotong University, Beijing, 100044, People's Republic of China\\
$^b$University of Primorska, Titov trg 4, 6000 Koper, Slovenia\\
$^c$University of Ljubljana, IMFM, Jadranska 19, 1000 Ljubljana, Slovenia}}
 \end{center}

\addtocounter{footnote}{-2}
\addtocounter{footnote}{1}  \footnotetext{Supported by KPCME
(106029) and SRFDP in China.}
\addtocounter{footnote}{1}  \footnotetext{Supported in part by
``Agencija za raziskovalno dejavnost Republike Slovenije'', research program P1-0285.

~*Corresponding author e-mail: ~dragan.marusic@guest.arnes.si}

\begin{abstract}

A  regular covering projection $\p\colon \tX \to X$ of connected  graphs 
is $G$-{\em admissible}  if $G$ lifts along $\p$. Denote by $\tG$ the lifted group,
and let $\CT(\p)$ be the group of covering transformations. The projection
is called  $G$-{\em split} whenever the extension $\CT(\p) \to \tG \to G$  splits. 
 In this paper,  split $2$-covers are considered, 
with a particular emphasis given to cubic symmetric graphs.
Supposing that $G$ is transitive on $X$,  a $G$-split cover is  said to be 
$G$-{\em split-transitive}
if  all complements $\bG \cong G$ of $\CT(\p)$ within $\tG$ are transitive
on $\tX$;  it is said to be $G$-{\em split-sectional} whenever for each complement $\bG$ 
there exists a $\bG$-invariant section of $\p$;  and it is called $G$-{\em split-mixed}
otherwise.

It is shown, when $G$ is an arc-transitive group, split-sectional and split-mixed  $2$-covers 
lead to canonical double covers. Split-transitive covers, however,
are considerably more difficult to analyze. 
For cubic symmetric graphs split $2$-cover are necessarily cannonical
double covers (that is, no 
$G$-split-transitive $2$-covers exist)
 when $G$ is  $1$-regular or $4$-regular.
In all other cases, that is, if $G$ is $s$-regular, $s=2,3$ or $5$,  
a necessary and sufficient condition for the existence of a transitive complement 
$\bG$ is given, 
and moreover, an infinite family  of  split-transitive $2$-covers
based on the alternating groups of the form $A_{12k+10}$ is constructed.

Finally, chains of consecutive $2$-covers, along which an arc-transitive group 
$G$ has successive lifts, are also considered. 
It is proved that in such a chain, at most two projections can be 
split. Further, it is shown that,  in the context of cubic symmetric graphs,
if exactly two of them are split, then one is split-transitive
and the other one is either split-sectional or split-mixed.

\end{abstract}

\bigskip
\begin{quotation}
\end{quotation}

%
%

\section{Introductory remarks}
\label{sec:intro}
\indent

 Let   $\p \colon \tX \to X$ be a  {\em regular  covering  projection} of connected (simple) graphs.
 Comparing symmetry properties of a given  base graph $X$ and a  covering graph $\tX$  has become
 quite an active area of research in recent years. The motivation stems from problems related to
 construction and classification of certain  classes of graphs and maps on surfaces,  counting the number
 of graphs in certain families,  producing lists of graphs with a given degree of symmetry,
 inductive approach to studying the structure of graphs via inspection of 
 smaller graphs arising as quotients  relative to a semiregular
 group of automorphisms etc. References are too numerous to be listed here, but see, for instance,
 \cite{AGS,CL89,CP96,XDK05,DMW98,FLP04,FK04,FK06,FKW05,FW03,GNS99,GLP04,H91,IP05,L01,L01a,LP00,LWX06,
 MMW04,M03,NS96,NS97,P93, RSJTW, ST97,WX06,Z73}.

 Questions related to symmetry properties of $X$ and $\tX$ are
 intimately linked with the problem of {\em  lifting and projecting automorphisms} \cite{BH72, D74, MMP, AM00}.
 Let  $g \in \Aut\,X$  and $\tilde g \in \Aut\,\tX$ be automorphisms  with the property that
 $ \p g = \tilde g \p$ (note that functions are composed on the right);
such an automorphism $\tilde g$ is known as a {\em lift} of $g$  while  $g$ is  the {\em projection}
of $\tilde g$ along $\p$.  If all elements
of a subgroup $G \leq \Aut\,X$
have a lift, then  the collection of all such lifts forms the {\em lifted group}  
$\tG \leq \Aut\, \tX$; the  projection $\p$ is then 
called {\em $G$-admissible}. In particular,
the group of {\em covering transformations}  $\CT(\p) \leq \Aut\,\tX$ is 
the lift of the trivial group, and  $\tG$ is isomorphic to an extension  of   $\CT(\p)$
by  $G$. We remark that any abstract group extension can be viewed as a lifting problem; 
this adds further motivation to  the topic. 

A $G$-admissible regular covering projection $\p \colon \tX \to X$
is called {\em split relative to $G$} (or  {\em $G$-split} for short)  whenever the extension
 $\CT(\p) \to \tG \to  G$  is split. This case deserves special attention: 
 since there exists
 a  complement $\bG \cong G$  of  $\CT(\p)$ within $\tG$  we can  compare actions  of two isomorphic groups,
  $G$ on $X$ and $\bG$ on $\tX$, where $\bG$ projects isomorphically onto $G$ along $\p$. A  restrictive
  situation such as this  makes it possible  to derive  a lot more information about  graphs and their symmetries. 
  
  However, the  analysis  can still be quite complicated.  This is due to the fact that  several complements 
   $\bG$ of  $\CT(\p)$  might
  exist within $\tG$, which moreover differ  in their actions  on $\tX$.
  Two particular cases of split covers seem to stand out.
  The first one  occurs  when there exists a complement $\bG$ which is {\em  sectional}.
  By this we mean that there is a {\em section} of $\tX$ -- 
  a set of vertices containing exactly one point from each fibre --
  which is invariant under the action of  $\bG$.
 The second  particular case occurs  when there exists a complement $\bG$ which is transitive on $\tX$ 
 ($G$ itself must then be transitive on $X$). Two particularly interesting extremal cases of split covers 
 can now be defined along these lines. 
 A projection $\p$ is called {\em split-sectional} if all  complements of $\CT(\p)$ within $\tG$ are sectional, and 
 is {\em split-transitive} if  all complements are transitive. A projection  is called {\em split-mixed} whenever
 both  sectional and transitive complements exist. Of course, it can also happen that  each complement of  $\CT(\p)$ within $\tG$
 is neither  sectional nor transitive. This fourth possibility is perhaps less interesting, and definitely the most difficult
 of all to analyze.

\medskip
 In this paper we restrict ourselves to {\em $2$-covers}, that is, to regular covering projections
 $\p\colon \tX \to X$   where 
 $\CT(\p)$ is isomorphic to  the cyclic group $\ZZ_2$, 
 with a particular emphasis given to such covering projections in the context of cubic symmetric graphs.
 Now,  if  the projection is  $G$-split, 
 the lifted group $\tG$ is necessarily isomorphic to the  direct product $G \times \ZZ_2$.
 In other words,  complements $\bG$ of $\CT(\p)$ are normal in $\tG$.
 In addition,  let us assume that $G$ is transitive on $X$.
Then   any  complement is either sectional or transitive in its action on $\tX$.
 Thus, apart form the possibility that $\p \colon \tX \to X$ is non-split relative to $G$, there are three possible {\em types}
 of $G$-split $2$-covers: split-sectional,  split-transitive, and split-mixed. 
 Obviously, these various kinds of lifts do depend on all  four items involved: the graph $X$, the cover $\tX$, the actual projection 
$\p \colon \tX \to X$, and the group $G$.
For example,  the cube $Q_3$ is obtained as a sectional  split $2$-cover  of $K_4$ relative to the 
alternating group $A_4$; on the other hand, $Q_3$ can be  viewed  as 
 a mixed  split $2$-cover  relative to the symmetric group $S_4$. For details see Section~\ref{sec-examples}
 where a  thorough analysis of all these possibilities is given 
in the context of some  cubic symmetric graphs of small order.

The rest of this article is organized as follows. In Section~\ref{sec:ex},
split lifts with a sectional a  complement  are considered  (for general graphs). 
In Section~\ref{subsec:transitive},  transitive complements are considered in the context of cubic symmetric graphs.
In Section~\ref{sec-height} we consider consecutive  $2$-lifts and show that within split covers only two such lifts are possible
for cubic symmetric graphs
We end  the paper with  Section~\ref{sec:remarks}, where we propose  a short list of problems  for future research.

%
%

\section{Examples}
\label{sec-examples}

\indent

\bigskip

We now give examples of all four types of $2$-covers, using 
the complete graph $K_4$ on $4$ vertices
and the Petersen graph $O_3$ as base graphs. Throughout this article 
most of our examples will be given in the context of  cubic symmetric graphs.
We will therefore  use  standard notation for these 
graphs from the extended Foster Census, see \cite{B88, CD02}, together with
commonly used names for some of  them where appropriate. 
By {\rm FnA}, {\rm FnB},  etc. we will refer to the corresponding graphs of order $n$ in
the Foster census of all cubic symmetric
graphs \cite{B88,CD02}, where the symbol {\rm FnA} is 
conveniently shortened to {\rm Fn} whenever a unique such graph exists.
For example, $K_4$ and $O_3$ are  the graphs $\F4$ and $\F10$, respectively.

Before starting with  examples, some additional terminology is in order.
A group of automorphisms  
$G$ of a graph $X$ is {\em arc-transitive} (also {\em symmetric}) 
if $G$ acts transitively
on the arcs in $X$,  and {\em $s$-regular} whenever  it acts regularly on the
$s$-arcs in $X$; here, an {\em $s$-arc} in $X$ is an ordered $(s+1)$-tuple
$(v_0,v_1,\ldots,v_{s-1},v_s)$ of vertices of $X$ with the property  that
$v_{i-1}$ is adjacent to $v_i$ for $1\leq i\leq s$ and
$v_{i-1}\not=v_{i+1}$ for $1\leq i<s$
(see \cite{Tu1}).
The graph $X$ is said to be
{\em $s$-regular} whenever  its full automorphism group $\Aut\ X$ acts
regularly on $s$-arcs in $X$. For instance, $K_4$ has only two 
arc-transitive groups of automorphisms,  the alternating group $A_4$ (which is $1$-regular) and the symmetric 
group $S_4$ (which is $2$-regular). Similarly,
the Petersen graph $O_3$ has only two arc-transitive groups of automorphisms,
$A_5$ and $S_5$, which are $2$- and $3$-regular, respectively.

Given a graph $X$, the {\em canonical double cover}    $\CDC(X) \to X$ is  a  $2$-cover   
which  can be reconstructed by 
the constant {\em  $\ZZ_2$-voltage assignment}  $\zeta(x) = 1$,    $x \in A(X)$. 
(See \cite{GT87} for an extensive treatment of 
voltage graphs.) Note that  $\CDC(X)$ is connected if and only if $X$ is connected and 
not bipartite. Note that 
 along the  canonical double cover any subgroup $G$
of the full automorphism group $\Aut\,X$ lifts to a group isomorphic to
$G \times \ZZ_2$, which  has a sectional complement $\bG \cong G$
to $\CT(\p)$, see Section~\ref{sec:ex}.

\medskip
\begin{example}
\label{ex:k4} {\rm The complete graph $K_4$ is not bipartite and
its canonical double cover $\p\colon CDC(K_4)\to
K_4$ is therefore connected. In fact,  $CDC(K_4)$   is isomorphic to the cube $Q_3$ of order $8$
($\F8$  in Foster notation); its  full automorphism group is  $S_4\times \ZZ_2$.

Clearly, $\p$ is a split $2$-cover relative to  the
$1$-regular group $A_4$ and relative to the $2$-regular group $S_4$ of $K_4$; 
both have sectional complements. Moreover,  no complement of $\p$ relative to $A_4$ is
transitive since  $8$ is not a divisor
of $|A_4|$.   Hence $\p$ is a  split-sectional $2$-cover relative to
$A_4$. 

Let $S_4\times\mz_2=S_4\times\langle a\rangle$ and let
$T=A_4\rtimes \langle ca\rangle$, where $S_4$ is a sectional
complement of $\p$ relative to the $2$-regular automorphism group
$S_4$ of $K_4$ and $c$ is an involution in $S_4 \setminus A_4$. Then, $T\cong
S_4$ is a transitive complement relative to $S_4$. This implies that
$\p$ is a  split-mixed  $2$-cover relative to $S_4$. }
\end{example}

\begin{example}
\label{ex:petersen}
{\rm
The Petersen graph $\F10$ is not bipartite and
hence it has a connected canonical double cover 
$\p_1\colon \CDC(\F10)\to \F10$. 
There are precisely two connected
cubic symmetric graphs of order $20$ (see \cite{B88,CN,CD02}), the dodecahedron
$\F20$A  and the Desargues graph $\F20$B.
The dodecahedron is  $2$-regular with the  full automorphism group $A_5\times\mz_2$.
The Deasargues graph, which is the canonical double cover $\CDC(\F10)$
is $3$-regular with the full automorphism group $S_5\times\mz_2$.
Note that the only arc-transitive subgroups of  $\Aut\, \F10$ are the
$2$-regular group $A_5$ and the $3$-regular group $S_5$. Clearly,
$A_5\times\mz_2$ contains only one subgroup isomorphic to $A_5$.
Thus, $\p_1$ is a  split-sectional  $2$-cover relative to $A_5$. 
On the other hand, a similar argument to the one given in  Example~\ref{ex:k4} 
shows that $\p_1$ is a  split-mixed $2$-cover relative to $S_5$.

The $2$-regular full automorphism group $A_5\times\mz_2$ of the
dodecahedron $\F20$A contains a normal subgroup $\mz_2$. The
quotient graph of  $\F20$A  corresponding to the orbits of $\mz_2$
must be the Petersen graph $\F10$ because there is only one connected
cubic symmetric graph of order $10$. Thus, $\F20$A  is a
split $2$-cover of $\F10$ relative to the
$2$-regular group $A_5$ of  $\F10$. 
Denote this $2$-cover by $\p_2\colon {\rm F20A}\to \F10$. 
The full automorphism group $S_5$ of $\F10$ cannot lift
along $\p_2$ because $\F20$A is $2$-regular. Since $\F20$A is
not bipartite, $A_5$ is transitive on $\F20$A, and the uniqueness
of $A_5$ in $A_5\times\mz_2$ implies that $\p_2$ is a transitive
$2$-cover relative to  $A_5$.

This example  needs some further comments which we here give without proof.
For details see Section~\ref{sec-height}. 
As $\F20$A  is not bipartite it  has the canonical double cover
 $\p_3\colon \F40\to {\rm F20A}$, 
where $\F40$ is the unique connected 
cubic symmetric graph of order $40$.  Note that $\F40$  is $3$-regular (see \cite{B88}) and
has $A_5\times\mz_2\times\mz_2$ as a $2$-regular group of automorphisms.
Then, $\p_3$ is a split-sectional  $2$-cover relative to the
$1$-regular automorphism group $A_5$  of {\rm F20A} and a 
split-mixed  $2$-cover of {\rm F20A} relative to the $2$-regular
automorphism group $A_5\times\mz_2$ of {\rm F20A}. 
Consider now the subgroup $A_5\times\mz_2\times\mz_2$ of
$\Aut\,\F40$.  One may choose a normal
subgroup $\ZZ_2$  in $A_5\times\mz_2\times\mz_2$ fixing the
bipartition  sets of  $\F40$ setwise.  The  quotient graph with respect to the action of this $\ZZ_2$ 
is bipartite,  and so it must be the Desargues graph ${\rm F20B}$.
Denote this projection, which is clearly not a CDC,  by $\p_4\colon \F40\to {\rm F20B}$. Then, $\p_4$
 is a  transitive $2$-cover of ${\rm F20B}$ relative to the $2$-regular automorphism
group $A_5\times\mz_2$ of ${\rm F20B}$. We obtain the following
figure: {\large
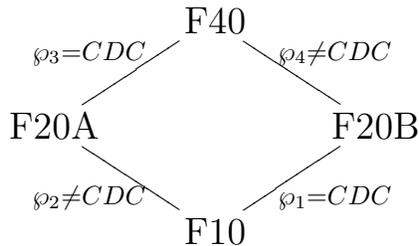
\begin{figure}[h!]
\centerline{\Large \xymatrix @-0pc {
*\txt{ } & *\txt{ }                   & \, \ar@{-}[dl]_{\p_3=CDC} \ar@{-}[dr]^{\p_4\ne CDC} \ar@{}[r]{\F40} & *\txt{ } \\
*\txt{ } & \, \ar@{}[l] {\rm F20A} \ar@{-}[dr]_{\p_2\ne CDC} & *\txt{ }                                          & \, \ar@{-}[dl]^{\p_1=CDC} \ar@{}[r] {\rm F20B} & *\txt{ }\\
*\txt{ } & *\txt{ }                   & \, \ar@{}[r]{\F10}                         & *\txt{ } \\
} } \caption{\label{fig:petersen} Consecutive split $2$-covers of
the Petersen graph relative to $A_5$.}
\end{figure}
}
Thus, there are two different `chains'  of consecutive $2$-covers
from the Petersen graph $\F10$ to the connected cubic symmetric 
graph $\F40$ relative to $A_5$. This idea will be used to prove that a `chain' of
consecutive split $2$-covers of a connected cubic graph relative to
an arc-transitive group of automorphisms  cannot have `length'  more than $2$
(see  Theorem~\ref{the:atmost2} in Section~\ref{sec-height}).
}
\end{example}

\begin{example}
\label{ex:infinitefam}
{\rm
In view of  Examples~\ref{ex:k4} and \ref{ex:petersen}  there exist 
transitive, sectional  and  mixed split $2$-covers. In
fact, infinitely many examples exist for each of these three cases.  By 
\cite[Theorem~4.4]{FK06}, there is an infinite family of
$2$-regular cubic graphs, denoted there by $EO_{p^3}$, where $p=\pm1 (\mod
5)$ is a prime. Each graph $EO_{p^3}$ is not bipartite and has $\mz_p^3\rtimes A_5$
as the full automorphism group. 
Since the cover is bipartite, any vertex-transitive subgroup of automorphisms
has a subgroup of index $2$ fixing the bipartition
sets setwise, and since the group  $\mz_p^3\rtimes A_5$
has no subgroup of index $2$,   the  canonical double cover of $EO_{p^3}$
 is a  split-sectional  $2$-cover relative to $\mz_p^3\rtimes A_5$. 

Again by \cite[Theorem~4.4]{FK06}, there is an infinite family of
$3$-regular cubic graphs arising as covers of the Petersen graph, denoted there by $EO_{p^6}$, where $p>5$ is a
prime. The graphs  $EO_{p^6}$ are not bipartite and have $\mz_p^6\rtimes S_5$
as the full automorphism group. 
Since $\mz_p^6\rtimes S_5$
has a transitive subgroup of index $2$, the canonical double cover of
$EO_{p^6}$  is a  split-mixed  $2$-cover relative to $\mz_p^6\rtimes S_5$
by Proposition~\ref{mixed}. Constructing infinitely many
 split-transitive  $2$-covers is a bit more complicated, see Theorem~\ref{transitive235} in Section~\ref{subsec:transitive}.
}
\end{example}

\begin{example}
\label{ex:nonsplit} 
{\rm 
To end this section we
would like to give an example of non-split $2$-cover, although 
this paper concentrates on  split $2$-covers.  

The M\"{o}bius-Kantor graph $\F16$ is a
$2$-cover of the $3$-cube $\F8$.  By \cite[Theorem]{FW03}, its  full
automorphism group $\Aut\,\F16 \cong \ZZ_2 \times S_4$ lifts. We
claim that this cover is non-split.  Suppose that $H=(\ZZ_2 \times S_4)\times \ZZ_2$
is a subgroup of $\Aut\,\F16$. Let $T$
be the unique non-trivial normal $2$-subgroup of $S_4$. Then  $T$
has order $4$,  and since $T$ is characteristic in $S_4$, it is normal
in $H$. It follows that the quotient graph of $\F16$
corresponding to the orbits of $T$ is the complete graph $K_4$,  and
$H/T\cong \ZZ_2\times S_3\times \ZZ_2$ can be viewed as a subgroup
of $\Aut\,K_4$  --  because $T$ is the kernel of the action of $H$ on the set
of orbits of $T$, which is impossible  since  $\Aut\,K_4 \cong S_4$.
This implies that the $2$-cover is non-split relative to the
automorphism group $\Aut\,\F8$. }
\end{example}

%
%

\section{Sectional complements in split $2$-covers}
\label{sec:ex}
\indent

Let  $\p \colon \tX \to X$  be a $2$-cover of connected graphs, and let $G \leq \Aut\,X$  be a vertex-transitive group that
lifts along $\p$ to $\tG$. The case when $\p$ is $G$-split with a sectional
complement can be conveniently
 described on the base graph  by an appropriate choice of voltages. In fact, the following
 holds (for a general version see \cite{AM00}).

 \begin{proposition}
 \label{nice}
 Let $\p \colon \tX \to X$ be a $2$-cover of connected graphs, and let $G \leq \Aut\,X$ be vertex-transitive.
 
{\rm  (1) }Then $G$ lifts along $\p$ if and only if, for any  voltage assignment by which the projection $\p$ is reconstructed,
 the set of all $0$-voltage closed walks is invariant under the action of $G$.

{\rm (2)} Moreover,
$\p$ is  $G$-split with a sectional
 complement  if and only if $\p$ can be reconstructed by a
voltage assignment $\zeta \colon A(X) \to \ZZ_2$   which satisfies  the following condition:
 for each automorphism $g \in G$ and each walk $W$ in $X$ we have
 $\zeta(W^g)=\zeta(W)$. 
  \end{proposition}

\begin{proof}
As for the first part, obviously  the condition that the set of $0$-voltage walks is invariant
under the action of $G$ is a necessary one. Similarly, it can also be seen  that it is  sufficient,
see \cite{AM00}.

We now prove the second part.
Suppose that $\p\colon \tX \to X$ is reconstructed by the voltage assignment in $\ZZ_2$
satisfying $\zeta(W^g)=\zeta(W)$ for all walks $W$ in $X$  and all $g \in G$. Then the set of
closed walks with  trivial voltage 
is clearly invariant under the action of $G$, and hence $G$ lifts. 
For each vertex $v$ of $X$, let $\fib(v) = \{v_0,v_1\}$ denote the fibre over $v$.
Choose a base vertex $b \in V(X)$. For each $g \in G$, let
$\bg$ be  the lift of $g$ which maps the vertex $b_0 \in \fib(b)$ to the vertex 
$b^g_0 \in \fib(b^g)$. Then $\bg$ preserves the set of vertices labelled by $0$, that is, the $0$-section.
Indeed, let $u \in V(X)$ be an arbitrary vertex, and let $W$ be the walk from $b$ to $u$ with $\zeta(W) = 0$.
(Note that such a walk always exists.)  Let $\bW$ be its lift with $b_0$ as the initial vertex. 
The terminal vertex of $\bW$ is $u_0 \in \fib(u)$  because $\zeta(W) = 0$.
The walk $W^g$ from $b^g$ to $u^g$ also has trivial voltage. Its lift 
$\bW^g$ starting at $b^g_0$ terminates at $u^g_0$. Since $\bg$ is the lift of $g$, it maps  the 
walk $\bW$ to the walk $\bW^g$. Hence $\bg$ maps $u_0$ to $u^g_0$. It follows that the $0$-section is invariant under 
the action of $\{\bg\ |\ g \in G\}$. Moreover  $\{\bg\ |\ g \in G\}$ must be a group, indeed a subgroup of index $2$ 
in $\tG$, as required.

To show the converse, let $\bG$ be a complement of $\CT(\p)$ preserving a section of $\p\colon \tX \to X$. 
This covering projection can be reconstructed by a voltage assignment in $\ZZ_2$ in such a way that this particular  section 
is labelled by $0$. Consider now an arbitrary walk $W$ in $X$, and let $W_0$ be its lift  in $\tX$ 
 starting at  the  corresponding vertex from the  $0$-section.   
 Since $\bG$ preserves the $0$-section and the $1$-section, it is obvious that
$W_0^{\bg}$ has the initial vertex in the $0$-section and its  terminal vertex belongs to the same section 
as the terminal vertex of $\bW$. Since the lift of $W^g$ starting in the $0$-section must be $W_0^{\bg}$, 
it follows that $\zeta(W^g) = \zeta(W)$, as required.
\end{proof}

\bigskip
A typical example illustrating Proposition~\ref{nice} is the canonical double cover 
$\CDC(X) \to X$. Recall that any group $G \leq \Aut\,X$ lifts to $CDC(X)$ 
as $G \times \ZZ_2$, and there exists  a sectional
complement. Further examples of split $2$-covers with a sectional complement 
which are not canonical double covers are given at the end of this section.
Note that a canonical double cover need not be split-sectional, for 
 transitive complements might exist  as well; hence it can be split-mixed, but clearly not split-transitive.
So let us  consider the question of when does  an arbitrary  $G$-split $2$-cover  with a sectional
complement have   another complement which is transitive on the
vertex set of the covering graph. The following holds.

\begin{proposition}
\label{mixed}
Let $\p\colon \tX \to X$ be a $G$-split $2$-cover
with a sectional complement. Then $\CT(\p)$ has a transitive complement within the  lifted group $\tG$
if and only if $G$ has a  vertex-transitive subgroup of index $2$.
\end{proposition}

\begin{proof}
Let $\tG = \CT(\p) \times \bG$, where $\bG$ has two orbits on the
vertex set of $\tX$.
Suppose that $\tG = \CT(\p) \times K$, where
$K$ acts transitively on the vertex set of $\tX$.  Since $\bG$ and $K$ are of index $2$ in $\tG$,
the intersection $K \cap \bG$ is of index $2$ in $\bG$ and in $K$. Consequently, the projection
$H \leq G$ of $\bG \cap K$ is a subgroup of index $2$ in $G$. Moreover, $H$ must be vertex-transitive.
Indeed, since the group $K \cap \bG$ is of index $2$ in $K$ it has at most two orbits on the
vertex set of $\tX$. In fact, as $K \cap \bG$ is contained in $\bG$, it cannot be transitive. So
$K \cap \bG$ must have two orbits and is therefore transitive on the set of all fibers. Consequently,
the  projection $H$ is transitive on $X$.

Conversely, let $H \leq G$ be a vertex-transitive subgroup of index $2$. Then there is an element
$g \in G$ such that $G = \langle H, g \rangle$. Denote by $c \in \CT(\p)$ the
nonidentity element. Then  $\tG = \langle c \rangle \times \bG$,  and the lift of
$H$ can be written as $\tH = \langle c \rangle \times \bH$ such that $\bH \leq \bG$.
Consider now the group $K = \langle \bH, c\bg \rangle$, where $\bg \in \bG$. Obviously,
$K$ is isomorphic to $G$ and acts transitively on the vertex set of $\tX$.
\end{proof}

\bigskip
See Example~\ref{ex:petersen} with  $F40 \to F20A$ and $G = A_5 \times \ZZ_2$ for an illustration  of the above proposition.
As already mentioned,  a  $G$-split $2$-cover  with a sectional complement
need not be a canonical double cover, see Construction~\ref{nonCDC} below. However, if $G$ is also edge-transitive 
(in addition to being  vertex-transitive),
then a $G$-split $2$-cover with a sectional complement is necessarily a canonical double cover.

\begin{proposition}
\label{niceCDC}
Let $\p\colon \tX \to X$ be a $G$-admissible $2$-cover of connected graphs.  If $G\leq \Aut\, X$ is  vertex and edge-transitive
then $\p$ is  $G$-split with a sectional complement if and only if $\p$ is the canonical double cover. (Equivalently,
if $\p$ is split, then it is   split-transitive relative to $G$ if and only if  $\p$  is not  the  canonical double cover).
\end{proposition}

\begin{proof}
If $\p$ is a canonical double cover, then clearly $\p$ is
$G$-split with a sectional complement. For the converse we may assume, by Proposition~\ref{nice},
that $\p$ is
reconstructed by a voltage assignment $\zeta : A(X) \to
\mathbb{Z}_{2}$ such that $\zeta(W^g)=\zeta(W)$ for all
walks $W$ and all $g \in G$. As $G$ is assumed
edge-transitive, all edges (and hence arcs) have equal voltage. Since 
 $\tX$ is assumed connected we have  $\zeta(x) = 1$ for all arcs, that
is,   $\p$ is the  canonical double cover.
\end{proof}

\bigskip
As long as we  choose to consider, say,  arc-transitive groups  (later on we shall in fact
restrict our considerations to arc-transitive cubic graphs), the question of
whether a given $G$-admissible $2$-cover is split with a sectional complement is solved
(in some sense). All we need is  to check  whether a  given $2$-cover is
indeed canonical. This  is algorithmically easy, as the
following proposition shows.

\begin{proposition}
\label{CDC}
Let $X$ be a connected graph and $\p\colon \tX  \to X$ a $2$-fold covering  projection arising from
a voltage assignment $\zeta \colon A(X) \to \ZZ_2$. Then $\p$ is the canonical double cover
if and only if each odd length cycle in
$X$ has voltage $1$ and each even cycle in $X$ has voltage $0$.
Moreover, it suffices to test this condition on the set of base cycles only.
\end{proposition}

\begin{proof}
It is easy to see that two voltage assignments $\zeta_1, \zeta_2 \colon A(X) \to \ZZ_2$ are equivalent
if and only if $\zeta_1(W) = \zeta_2(W)$ for each closed walk $W$, and the proof is immediate.
\end{proof}

\bigskip

We now turn to the question of existence of split $2$-covers of a vertex-transitive graph $X$ with sectional complement which are not 
canonical double covers. In view of Proposition~\ref{niceCDC}, this can only happen if the vertex-transitive group $G$ in question is 
not edge-transitive.  We  give below two examples of split $2$-covers (of cubic vertex-transitive graphs)
which are not canonical double covers.

\begin{figure}[h!]
\begin{minipage}[b]{0.4\linewidth} 
\centering
\includegraphics[width=4cm]{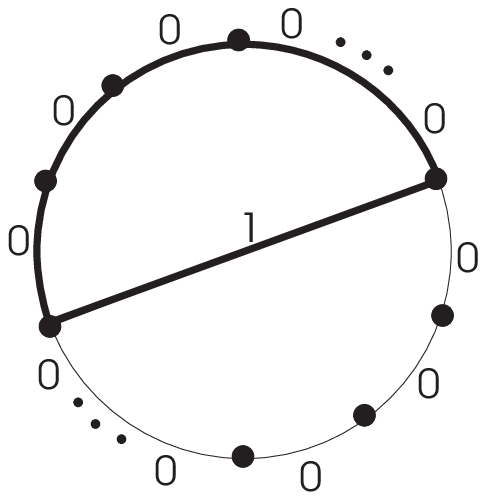}
\caption{\label{fig:2}}
\end{minipage}
\hspace{0.08cm} 
\begin{minipage}[b]{0.6\linewidth}
\centering
\includegraphics[width=6cm]{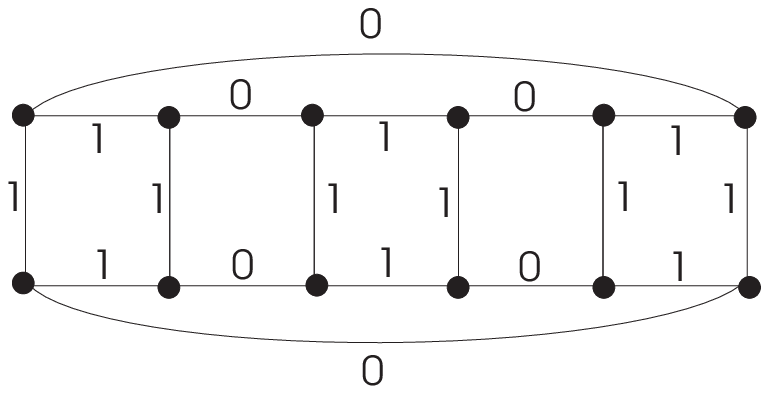}
\caption{\label{fig:3}  }
\end{minipage}
\end{figure}

\begin{construction}
\label{nonCDC}
{\rm
Let $n$ be odd and let $X$ be the circulant
$X \cong \Cay(\ZZ_{2n}, \{1,2n-1,n\}\})$.
Define the voltage assignment
$\zeta : A(X)\rightarrow \ZZ_2$  in such a way that all edges $[i,i+1]$ $i \in\ZZ_{2n}$ receive voltage $0$
and edges  $[i,i+n]$ receive voltage $1$.
Moreover, by Proposition~\ref{CDC}, the respective $2$-cover $\p\colon \tX \to X$ is not the  canonical double cover
 since the cycle $(0,1,2\dots n,0)$ 
of $X$ has even length but voltage $1$
(see Figure~\ref{fig:2}).
Let $G =\Aut\, X \cong D_{4n}$, let $H \cong \ZZ_{2n}$ be the  cyclic subgroup of index $2$ in $G$, 
and let $K \cong D_{2n}$ be the dihedral subgroup of index $2$ in $G$.   
Obviously, the group $G$ lifts by Proposition~\ref{nice}, and the corresponding $2$-cover is split. Moreover,
by Proposition~\ref{niceCDC}, the projection $\p$ is split-sectional relative to $H$ and relative to $K$, 
and split-mixed relative to $G$. 
}
\end{construction}

\begin{construction}
\label{nonCDC-another}
{\rm
Let $n$ be odd, and let
$X=C_{2n}\square\, K_2$ be the cartesian
product of $C_{2n}$ with  $K_2$,  viewed as the 
Cayley graph $\Cay(G, \{b, ab, c\})$, where the group
$G= \langle a, b, c \ | \ a^n=b^2=c^2=1, a^b=a^{-1}, a^c=a, b^c =b \rangle$
is isomorphic to $D_{2n}\times \ZZ_2$.
Define the voltage assignment 
$\zeta : A(X) \rightarrow \ZZ_2$  'consistent' with the regular action of $G$
in such a way that all $b$-edges and all $c$-edges receive voltage $1$ and all 
$ab$-edges receive voltage $0$.
The corresponding $2$-cover is not the  canonical double cover because $X$ 
contains an even length cycle with voltage $1$. 
In fact,  any of the two $2n$-cycles obtained from  $b$-edges and $ab$-edges is of this kind
(see Figure~\ref{fig:3}).
Also, this $2$-cover has a sectional  complement relative to $G$.
It is split-sectional.
}
\end{construction}

%
%

\bigskip
\section{Split $2$-covers with a transitive complement}
\label{subsec:transitive}

We start this section by constructing a graph from a finite group
$G$ relative to a subgroup $H$ of $G$ and a union $D$ of some double
cosets of $H$ in $G$ such that $D^{-1}=D$, the so called coset
graph \cite{S, Mi}.

The {\em coset graph} $\Cos(G,H,D)$ of $G$ with respect to
$H$ and $D$ is defined to have vertex set $[G:H]$, the set of right
cosets of $H$ in $G$, and edge set $\{(Hg,Hdg)\ |\ d\in D\}$. The
graph has valency $|D|/|H|$ and is connected if and only if $D$
generates the group $G$. For $g\in G$, denote by $R(g)_H$ the right
multiplication of $g$ on $[G:H]$, that is, $R(g)_H:$ $Hx\mapsto
Hxg$, $x\in G$. Clearly, $R(g)_H$ is an automorphism of
$\Cos(G,H,D)$. Let $R(G)_H=\{R(g)_H\ |\ g\in G\}$. Then $R(G)_H\leq
\Aut(\Cos(G,H,D))$. Note that it may happen that $R(g)_H$ is the
identity automorphism of $\Cos(G,H,D)$ for some $g\not=1$ in $G$.
Actually, $R(G)_H\cong G$ if and only if $H_G=1$, where $H_G$ is the
largest normal subgroup of $G$ in $H$. Since $R(G)_H\leq
\Aut(\Cos(G,H,D))$, we have that $\Cos(G,H,D)$ is vertex-transitive and $R(G)_H$
acts arc-transitively  on the coset graph if and only if $D$ is a
single double coset. (Note that the concept of a coset graph is equivalent
to the concept of an orbital graph \cite{S67, P90}.)

Conversely, let $X$ be a graph and let $A$ be a vertex-transitive subgroup
of $\Aut(X)$. By \cite{S}, the graph $X$ is isomorphic to a coset graph
$\Cos(A,H,D)$, where $H=A_u$ is the stabilizer of $u\in V(X)$ and
$D$ consists of all elements of $A$ which map  $u$ to its neighbors.
It is easy to show that $H_A=1$ and that $D$ is a union of some double
cosets of $H$ in $A$  satisfying  $D=D^{-1}$. Assume that $A$ acts
arc-transitively on $V(X)$ and that $g\in A$ interchanges $u$ and one of
its neighbors. Then, $g^2\in H$ and $D=HgH$. The valency of $X$ is
$|D|/|H|=|H:H\cap H^g|$.

 Let $X$ be a connected cubic symmetric  graph with $G$ as an
$s$-regular subgroup of $\Aut(X)$. By
\cite[Proposition~2-Proposition~5]{DM}, the stabilizer $G_v$ of
$v\in V(X)$ in $G$ is isomorphic to $\mz_3$, $S_3$,
$S_3\times\mz_2$, $S_4$, or $S_4\times \mz_2$  for $s=1,2,3,4$ or
$5$, respectively.

Now we consider split $2$-covers of a connected cubic graph relative
to an $s$-regular  group  of automorphisms with a transitive
complement. Let $X$ have a connected split $2$-cover $Y$ such that
an $s$-regular  group of automorphisms $G$ lifts to $\mz_2\times G$ with
$G$ acting  on $Y$ transitively. Clearly, $|V(Y)|=2|V(X)|$. If $s=1$
then $|G|=3|V(X)|$ and since $G$ is transitive on $Y$, $|Y|$ is a
divisor of $|G|$, that is, $2|V(X)| \mid 3|V(X)|$, which is
impossible. If $s=4$ then $\mz_2\times G$ is $4$-regular and $G$ is
$3$-regular, which is also impossible because a $3$-regular
group of automorphisms  cannot be contained a $4$-regular group of automorphisms
 (see \cite{DM}). This implies the following theorem.

\begin{theorem}\label{transitive14} Let $X$ be a connected cubic graph
with $G$ as an $s$-regular automorphism group,  where $s=1,4$. Then
each $G$-split $2$-cover must be the canonical double cover.
\end{theorem}

\begin{theorem}\label{transitive235} Let $X$ be a connected cubic graph
with $\{u,v\}$ as an edge and with $G$ as an $s$-regular
group of automorphisms,  where $s=2,3,5$, and let $H=G_u$ be the
stabilizer of $u\in V(X)$ in $G$. Then, $X$ has a connected $G$-split
$2$-cover with a transitive complement if and only
if there is an element  $b\in G$ interchanging $u$ and $v$ and $H$ has a
subgroup $L$ of index $2$ such that
\begin{eqnarray}
b^2\in L,
\quad
\langle b,L\rangle=G,\quad \hbox{{\rm  and}}\quad 
|L:L^b\cap L|=3.  \label{cond}
\end{eqnarray}
Moreover, the covering graph is isomorphic to the coset graph $\Cos(G,L,LbL)$.
\end{theorem}

\begin{proof} 
Suppose first that $X$ has a connected split
$2$-cover $Y$ such that $G$ lifts to $\bG \times \ZZ_2$ with $\bG$ acting
on $Y$ transitively. In the sequel we may identify $\bG$ with $G$ for convenience.
Then, $G$ is $(s-1)$-regular on $V(Y)$ because
it is $s$-regular on $V(X)$. Since $s=2,3,5$, the group $G$ is arc-transitive
on $Y$. Let $\{u_1,u_2\}$ and $\{v_1,v_2\}$ be the fibres of $u$ and
$v$ in the $2$-cover $Y$ of $X$. By identifying each vertex in
$V(X)$ with its fibre in $Y$, one has $H=G_u=G_{\{u_1,u_2\}}$, where
$G_{\{u_1,u_2\}}$ is the subgroup of $G$ fixing the fibre
$\{u_1,u_2\}$ setwise. Since $u$ and $v$ are adjacent, $u_1$ is
adjacent to a vertex in $\{v_1,v_2\}$, say $v_1$. The
arc-transitivity of $G$ on $Y$ implies that there is an element $b\in G$
interchanging $u_1$ and $v_1$. It follows that $b$ interchanges $u$
and $v$. Let $L=G_{u_1}$. Since $Y$ is a connected cubic graph, one
has $G=\langle b,L\rangle$, $b^2\in L$ and $|L:L^b\cap L|=3$.
This shows that (\ref{cond}) is  a necessary condition.

For the converse, suppose that 
there is an element  $b\in G$ interchanging  $u$ and $v$ and that $H$ has
a subgroup $L$ of index $2$ such that (\ref{cond}) is satisfied.

Since $G$ is arc-transitive on
a given connected cubic graph $X$, one has $H_G=1$, $|H:H^b\cap H|=3$, 
and $X\cong\Cos(G,H,HbH)$.  Without loss of generality we may identify $X$ with 
this coset graph, 
$X=\Cos(G,H,HbH)$.  Set $\tX =\Cos(G,L,LbL)$. Clearly, $\tX$  is a  connected cubic
graph with $G$ as an $(s-1)$-regular automorphism subgroup since
$H_G=1$ implies that $L_G=1$. Thus, $L$ is isomorphic to $\mz_3, S_3, S_4$ for
$s=2,3,5$, respectively. Since $H$ is isomorphic to $S_3, S_3\times \mz_2,
S_4\times \mz_2$ for $s=2, 3, 5$, respectively, there is an involution $r \in H$
such that $H=L\rtimes \langle r\rangle$, where $r$ normalizes $L$. 
Note that for $s=3,5$ the involution $r$ centralizes $L$ (and hence $H = L \times \la r\ra$),
but for $s=2$ it does not.
Set
\begin{eqnarray}
K= G \times \langle c\rangle,\ \ D=\langle L,rc\rangle,  \label{cond2}
\end{eqnarray}
where $\la c \ra \cong \ZZ_2$.
Then, $D\leq K$, $\langle DbD\rangle=\langle D,b\rangle=K$ and
$(DbD)^{-1}=DbD$ because $\langle L,b\rangle=G$ and $b^2\in L$.
Since $rc$ normalizes $L$ (but does not centralize) for $s=2$,  and
centralizes  $L$ for $s=3,5$, one has $D\cong S_3, S_3\times\mz_2,
S_4\times \mz_2$ for $s=2,3,5$, respectively. We claim that 
$$%
\overline{X}:=\Cos(K,D,DbD)
$$%
is a connected cubic graph with $K$ as
an  $s$-regular group of automorphisms. To prove this  we need  to show that
$D_K=1$ and that $|D:D^b\cap D|=3$.

Suppose on the contrary that $D_K\not=1$. Then, $D_K\lhd D$, $D_K\cap G\lhd G$,  and
$D_K\cap L\lhd L$. Let $s=2$. Then  $D\cong S_3$,  and since $D_K\lhd
D$, one has that $D_K=L$  is the unique Sylow $3$-subgroup of $D$. Hence 
$D_K\leq H$, which contradict the fact that $H_G=1$. Let $s=3$. 
Assume that
$N\leq D_K$ is a minimal normal subgroup of $K$. Then $N$ is an
elementary abelian $2$- or $3$-group. Since $D=L\times\langle rc\rangle\cong S_3\times \mz_2$, 
the group $N$ is either the Sylow $3$-subgroup of
$L$,  or else $N=\langle rc\rangle$.  In  the former case  one has $N\leq H_G$, a
contradiction. As for  the latter case,  we have $rc\in  Z(K)$ 
because $|N|=2$ and hence $\langle r\rangle\leq Z(G)$. It
follows that $\langle r\rangle \leq H_G$, a contradiction. Let $s=5$.
Then $D=L\times \langle rc\rangle\cong S_4\times \mz_2$. Clearly,
$L\cong S_4$ has only one non-trivial normal subgroup of order $4$, say $T \cong \ZZ_2 \times \ZZ_2$.
Let $M\leq D_K$ be a minimal normal subgroup of
$K$. Then  $M$ is a $2$- or a $3$-group. Since $M\lhd D$,  the group $M$ cannot
be a $3$-group, and one may easily obtain that $M=T$, $ M = T\times \langle rc\rangle$, or $M =\langle rc\rangle$. 
For the first two cases,
$T=M\cap G\lhd G$,  and for the last case $\langle r\rangle \lhd G$,
contrary to $H_G=1$. Thus, $D_K=1$.

To  prove $|D:D^b\cap D|=3$, it suffices to show  that $D^b\cap D$ is a
Sylow $2$-subgroup of $D$. Since $|H:H^b\cap H|=|L:L^b\cap L|=3$, the groups 
$H^b\cap H$ and $L^b\cap L$ are Sylow $2$-subgroups of $H$ and $L$, respectively.
Since $b^2\in L\leq H$, one has $(H^b\cap H)^b=H^b\cap H$ and
$(L^b\cap L)^b=L^b\cap L$, that is, $b$ normalizes  $H^b\cap H$ and
$L^b\cap L$. Let $s=2$. Then $L\cong\mz_3$, $H=L\rtimes\langle r\rangle$,  
and $D=L\rtimes \langle rc\rangle\cong S_3$. This implies
that $L^b\cap L=1$ and $H^b\cap H\cong \mz_2$. Assume that $H^b\cap
H=\langle h\rangle$. Since $b$ normalizes  $H^b\cap H$ and since $c\in Z(K)$,
one has $h^b=h$ and $(ch)^b=ch$. Note that $h$ is a product
of $r$ and an element in $L$; this holds since  $H\cong S_3$. Thus, $ch\in D$, 
and so $D^b\cap D$ contains the Sylow $2$-subgroup $\langle ch\rangle$ of $D$. 
Suppose $D^b\cap D=D$. Then $D^b=D$,  and hence $b$
normalizes  the unique Sylow $3$-subgroup $L$ of $D$, that is,
$L^b=L$. Thus, $L\leq H^b\cap H$, which is impossible because
$H^b\cap H$ is a $2$-group. It follows that $D^b\cap D$ is a Sylow
$2$-subgroup of $D$, as required. So let us consider the remaining two cases where  $s=3$ or $5$. 
Then   $L\cong S_3$,
$H=L\times \langle r\rangle$ and $D=L\times \langle rc\rangle$ for
$s=3$, and $L\cong S_4$, $H=L\times \langle r\rangle$ and $D=L\times
\langle rc\rangle$ for $s=5$. In both cases,  let $T$ be a Sylow
$2$-subgroup of $L$. Then, $H$ has a unique Sylow $2$-subgroup
containing $T$, that is, $T\times \langle r\rangle$. Since $L^b\cap L\leq H^b\cap H$ 
and the groups $H^b\cap H$ and  $L^b\cap L$ are Sylow
$2$-subgroups of $H$ and $L$, respectively, one has $H^b\cap H=(L^b\cap L)\times
\langle r\rangle$. Thus, $(L^b\cap L)\times \langle r\rangle=H^b\cap
H=(H^b\cap H)^b=((L^b\cap L)\times \langle r\rangle)^b=(L^b\cap
L)\times \langle r^b\rangle$. It follows that there exists an
$\ell\in L^b\cap L$ such that $r^b=\ell r$. So, $(L^b\cap L)\times
\langle rc\rangle)^b=(L^b\cap L) \times \langle r^bc\rangle=(L^b\cap
L)\times \langle \ell rc\rangle=(L^b\cap L)\times\langle rc\rangle$,
implying that $D^b\cap D$ contains the Sylow $2$-subgroup $(L^b\cap
L)\times\langle rc\rangle$ of $D$. Note that a Sylow $3$-subgroup of
$D$ is also a Sylow $3$-subgroup of $H$. This implies that $D^b\cap
D\not=D$ because $H^b\cap H$ is a Sylow $2$-subgroup of $H$. It
follows that $D^b\cap D$ is a Sylow $2$-subgroup of $D$, as required.

Thus, the claim is true, that is, $\overline{X}=\Cos(K,D,DbD)$ is a
connected cubic graph with $K$ as an  $s$-regular group of automorphisms.
Recall that 
$K=G\times\langle c \rangle$, $H=L\rtimes \langle r\rangle$,  and $D=L\rtimes\langle rc \rangle$.
For $g\in G$ we have
$Dgc=Drcgc=Drg$. Thus, $G$ is transitive on $\overline{X}$. Let
$\underline{X}$ be the quotient graph of $\overline{X}$
corresponding to the orbits of $\langle c\rangle$. Then 
$\overline{X}$ is a $2$-cover of $\underline{X}$,  with $G$ projecting
to a group isomorphic to $G$. Recall that  $X=\Cos(G,H,HbH)$. For a
vertex $Dk\in V(\overline{X})$, denote by $\underline{Dk}$ the vertex
of \underline{X} corresponding to $Dk$, that is,
$\underline{Dk}=\{Dk,Dkc\}$, the orbit of $\langle c\rangle$
containing $Dk$. Define a map from $X$ to \underline{X} by 
$\a$: $Hg\mapsto \underline{Dg}$ for $g\in G$. Assume $Hg_1=Hg_2$. Since
$H=L\cup Lr$, one has $Lg_1\cup Lrg_1=Lg_2\cup Lrg_2$. It follows
that either $Lg_1=Lg_2$ and $Lrg_1=Lrg_2$ or $Lg_1=Lrg_2$ and
$Lrg_1=Lg_2$. This implies that 
$\underline{Dg_1}=\{Dg_1,Dg_1c\}=
    \{Lg_1\cup Lrg_1c, Lrg_1\cup Lg_1c\}=\{ Lg_2\cup Lrg_2c, Lrg_2\cup Lg_2c\}=\underline{Dg_2}$ 
because $D=L\cup Lrc$. Thus, $\a$ is well
defined. By transitivity of $G$ on $\overline{X}$, the mapping  $\a$ is
surjective and so bijective because $|V(X)|=|V(\underline{X})|$.
Take an edge $(Hg,Htg)$ in $X$,  where $g\in G$ and $t\in HbH$. Since
$H=\langle L, r\rangle$ and $D=\langle L, rc\rangle$, one has $t\in
DbD$ or $tc\in DbD$. If $t\in DbD$ then $\{Dg, Dtg\}$ is an edge of
$\overline{X}$, and if $tc\in DbD$ then $\{Dg, Dtgc\}$ is an edge. In
both cases, $(\underline{Dg}, \underline{Dtg})$ is an edge of
\underline{X}. Thus, $\a$ is an isomorphism from $X$ to
\underline{X}. This means that $\overline{X}$ is a $2$-cover of $X$, 
with $G$ projecting to a group isomorphic to $G$.

As  with  the isomorphism  $\a$, one can  easily show that the map from $\Cos(G,L,LbL)$ to
$\overline{X}=\Cos(K,D,DbD)$, defined by 
$\b$: $Lg\mapsto Dg$,   is an isomorphism. This completes the
proof. \hfill
\end{proof}

\vskip 0.3cm

We say that an ordered  pair groups $(H,K)$ is isomorphic to an ordered pair groups
$(H',K')$ if $H\cong H'$ and $K\cong K'$; this is  denoted by
$(H,K)\cong (H',K')$ for short. By the constructions of the automorphisms  $\a$
and $\b$ in Theorem~\ref{transitive235}, we have

\begin{corollary}\ \ \label{corollary} Let $L\leq H\leq G$ with $H_G=1$. Suppose
that $(H,L)\cong
(S_3,\mz_3), (S_3\times \mz_2,S_3)$ or $(S_4\times\mz_2,S_4)$. Then,
$G$ has an element $b$ satisfying $b^2\in L$, $\langle b,L\rangle=G$,
 and $|H:H^b\cap H|=|L:L^b\cap L|=3$   if and
only if the projection $\p: Lg\mapsto Hg$ from the coset graph
$\Cos(G,L,LbL)$ to the coset graph $\Cos(G,H,HbH)$ is a split
$2$-cover with $G$ as a transitive complement.
\end{corollary}

In the remainder of this section  we construct an infinite family of
split-transitive $2$-covers of  cubic symmetric graphs  relative to $2$-regular
groups of automorphisms. We will need two classical results regarding transitive permutation groups, 
by Jordan and Marggraf \cite{WI}. Let $G$ be
a transitive permutation group on a set $\O$. A nonempty subset $\D$
of $\O$ is called a {\em block} for $G$ if,  for each $g\in G$,  either
$\D^g=\D$ or $\D^g\cap\D=\phi$. A group $G$ is said to be {\it primitive} 
if $G$ has no block $B$ such that $1<|B|<|\O|$.

\begin{proposition}
{\rm\cite[Thm.~13.8]{WI}} \label{Marggraf} A primitive group of
degree $n$, which contains a cycle of degree $m$ with $1<m<n$, is
$(n-m+1)$-fold transitive.
\end{proposition}

\begin{proposition}{\rm\cite[Thm.~13.9]{WI}}\ \ \label{Jordan}
Let $p$ be a prime and $G $ a primitive group of degree $n=p+k$ with
$k\geq 3$. If $G $ contains an element of degree $p$ and order $p$,
then $G $ is either alternating or symmetric.
\end{proposition}

The construction of split-transitive $2$-covers is based on the alternating group $A_{2k+10}$, 
where $k$ is a nonnegative integer. Let us first define the following permutations in $A_{12k+10}$:
$$%
\begin{array}{rl}
a_k=&\displaystyle\prod _{i=1}^{4k+3}(3i-2\ 3i-1\ 3i), \, \, \\
r_k=&\left(\displaystyle\prod _{i=1}^{2k+1}(6i-2\ 6i-1)(6i+1\
6i+2)\right)(1\ 2), \, \, \\
b_k=&\left(\displaystyle\prod _{i=1}^{2k+1}(6i-3\ 6i)(6i-2\
6i+1)(6i-1\ 6i+2)\right)(12k+9\ 12k+10).
\end{array}
$$%
Clearly,
$L_k=\langle a_k\rangle\cong\mz_3$,
$H_k=\langle a_k,r_k\rangle\cong S_3$,  and
$H_k^{b_k}\cap H_k=\langle r_k\rangle\cong \mz_2$.

\begin{theorem} \ \ \label{transitive} The projection $\p: L_kg\mapsto H_kg$, $g\in
 A_{12k+10}$, is a split-transitive  $2$-cover from the
coset graph $\Cos(A_{12k+10}, L_k, L_kb_kL_k)$ to the coset graph
$\Cos(A_{12k+10}, H_k, H_kb_kH_k)$ relative to $A_{12k+10}$. 
\end{theorem}

\begin{proof} If $A_{12k+10}=\langle a_k,b_k\rangle$ then the theorem is true by
Corollary~\ref{corollary} because the simplicity of $A_{12k+10}$
implies that there is a unique complement in the $2$-cover $\p$. Let
$G_k=\langle a_k,b_k\rangle$. To finish the proof, it suffices to
show that $G_k=A_{12k+10}$.

Let $i_1, i_2, \ldots,  i_\ell, j_1,  j_2, \ldots,  j_m$ be distinct
numbers. Let $x=(i_1\ i_2\ \ldots\ i_\ell)$ and $y=(j_1\ j_2\ \ldots\ j_m)$ 
be cycle permutations with the first entry $i_1$ and
$j_1$ distinguished, respectively. By $(i_1\ i_2\ \ldots\ i_\ell\ y)$ or $(x\ y)$ we denote the cyclic
permutation ($i_1\ i_2\ \ldots\ i_\ell\ j_1\ j_2\ \ldots\ j_m$).

It is easy to see that $G_k$ is transitive on
$\O=\{1,2,\cdots,12k+10\}$. Now, use an induction on $k$ to claim
that
\begin{equation}
\begin{array}{l}
b_kb_k^{a_k}=(e_k\ f_k)\end{array}  \label{ef}
\end{equation}
where $e_k$ is a cycle of length $4k+4$ on $\{3\ell+1\ |\ 0\leq
\ell\leq 4k+3\}$ with the first entry $12k+7$ and the last entry
$12k+10$ and $f_k$ is a cycle of length $4k+3$ on $\{3\ell\ |\ 1\leq
\ell\leq 4k+3\}$ with the first entry $12k+6$ and the last entry
$12k+9$.

If $k=0$ then $b_0b_0^{a_0}=(7\ 1\ 4\ 10\ 6\ 3\ 9)=(e_0\ f_0)$,
where $e_0=(7\ 1\ 4\ 10)$ is a cycle on $\{1,4,7, 10\}$ and $f_0=(6\
3\ 9)$ is a cycle on $\{3,6,9\}$. The claim is true. By induction
hypothesis, assume that (\ref{ef}) is true for $k=m\geq 1$, that is,
$b_mb_m^{a_m}=(e_m\ f_m)$,  where $e_m$ is a cycle of length
$4m+4$ on $\{3\ell+1\ |\ 0\leq \ell\leq 4m+3\}$ with the first entry
$12m+7$ and the last entry $12m+10$ and $f_m$ is a cycle of length
$4m+3$ on $\{3\ell\ |\ 1\leq \ell\leq 4m+3\}$ with the first entry
$12m+6$ and the last entry $12m+9$.

Let $k=m+1$. Set
$$
\begin{array}{ll}
a'=  &(12m+10\ 12m+11\ 12m+12),\\
a''= &(12m+13\ 12m+14\ 12m+15)(12m+16\ 12m+17\ 12m+18)\\
     & (12m+19\ 12m+20\ 12m+21),\\
b'= &(12m+9\ 12m+13\ 12m+10\ 12m+12),\\
b''=&(12m+11\ 12m+14)(12m+15\ 12m+18)\\
 & (12m+16\ 12m+19)(12m+17\ 12m+20)(12m+21\ 12m+22).
\end{array}
$$
Then, $a_{m+1}=a_ma'a''$ and $b_{m+1}=b_mb'b''$. Furthermore,
$a_ma''=a''a_m$, $a_mb''=b''a_m$, $b_mb''=b''b_m$ and
$b_ma''=a''b_m$. Thus,
$b_{m+1}a_{m+1}=b_mb'b''a_ma'a''=b_ma_m(b')^{a_m}b''a'a''=b_ma_mc_1$,
where $c_1=(b')^{a_m}b''a'a''$. Note that $(b')^{a_m}=(12m+7\
12m+13\ 12m+10\ 12m+12)$ and by computation, one has
$$
\begin{array}{ll}
c_1=&(12m+7\ 12m+14\ 12m+12)(12m+17\ 12m+21\ 12m+22\ 12m+19)\\
&(12m+11\ 12m+15\ 12m+16\ 12m+20\ 12m+18\ 12m+13).
\end{array}
$$
Similarly,
$a_{m+1}b_{m+1}a_{m+1}=a_mb_ma_m(a')^{b_ma_m}a''c_1=a_mb_ma_mc_2$
where
$$
\begin{array}{ll}
c_2=&(12m+7\ 12m+15\ 12m+11)(12m+18\ 12m+20\ 12m+22\ 12m+19)\\
&(12m+12\ 12m+14\ 12m+16\ 12m+21\ 12m+17\ 12m+13);
\end{array}
$$
$a_{m+1}^2b_{m+1}a_{m+1}=a_m^2b_ma_m(a')^{a_mb_ma_m}a''c_2=a_m^2b_ma_mc_3$
where
$$
\begin{array}{ll}
c_3=&(12m+11\ 12m+14)(12m+12\ 12m+15)(12m+13\ 12m+16)\\
&(12m+17\ 12m+20)(12m+18\ 12m+21)(12m+19\ 12m+22);
\end{array}
$$
$b_{m+1}a_{m+1}^2b_{m+1}a_{m+1}=b_ma_m^2b_ma_m(b')^{a_m^2b_ma_m}b''c_3=b_ma_m^2b_ma_mc_4$
where
$$
\begin{array}{ll}
c_4=&(12m+6\ 12m+16\ 12m+22\ 12m+18\ 12m+12)\\
&(12m+7\ 12m+15\ 12m+21\ 12m+19\ 12m+13).
\end{array}
$$
Thus, $b_{m+1}b_{m+1}^{a_{m+1}}=b_mb_m^{a_m}c_4=(e_m\ f_m)c_4$. 
Recall 
that $e_m$ is a cycle of length $4m+4$ on $\{3\ell+1\ |\ 0\leq
\ell\leq 4m+3\}$ with the first entry $12m+7$ and the last entry
$12m+10$ and $f_m$ is a cycle of length $4m+3$ on $\{3\ell\ |\ 1\leq
\ell\leq 4m+3\}$ with the first entry $12m+6$ and the last entry
$12m+9$. It follows that $ b_{m+1}b_{m+1}^{a_{m+1}}=(12m+19\ 12m+13\
e_m\ 12m+16\ 12m+22\ 12m+18\ 12m+12\ f_m\ 12m+15\ 12m+21)$. Set
$e_{m+1}=(12m+19\ 12m+13\ e_m\ 12m+16\ 12m+22)$ and
$f_{m+1}=(12m+18\ 12m+12\ f_m\ 12m+15\ 12m+21)$. Then,
$b_{m+1}b_{m+1}^{a_{m+1}}=(e_{m+1}\ f_{m+1})$, where $e_{m+1}$ is a
cycle of length $4(m+1)+4=4m+8$ on $\{3\ell+1\ |\ 0\leq \ell\leq
4(m+1)+3\}$ with the first entry $12(m+1)+7=12m+19$ and the last
entry $12(m+1)+10=12m+22$ and $f_{m+1}$ is a cycle of length
$4(m+1)+3=4m+7$ on $\{3\ell\ |\ 1\leq \ell\leq 4(m+1)+3\}$ with the
first entry $12(m+1)+6=12m+18$ and the last entry
$12(m+1)+9=12m+21$. Thus, (\ref{ef}) is true for any $k\geq 0$.

Recall  that $\O=\{1,2,\ldots, 12k+10\}$. Let $\O_1=\{3\ell+1\ |\ 0\leq
\ell\leq 4k+3\}$, $\O_2=\{3\ell+2\ |\ 0\leq \ell\leq 4k+2\}$ and
$\O_3=\{3\ell\ |\ 1\leq \ell\leq 4k+3\}$. By (\ref{ef}),
$b_kb_k^{a_k}=(e_k\ f_k)$ is a cyclic permutation on $\O_1\cup\O_3$
of length $8k+7$ with $e_k$ and $f_k$ cyclic permutations on $\O_1$
and $\O_3$ respectively, and it fixes the set $\O_2$ pointwise.
Thus, $(b_kb_k^{a_k})^{a_k}=(e'_k,f'_k)$ is a cyclic permutation of
length $8k+7$ on $\O_1\cup\O_2$ with $e'_k$ and $f'_k$ cyclic
permutations on $\O_2$ and $\O_1$ respectively. Remember that $G_k $
is transitive on $\O$. We now prove that $G_k$ is primitive. Let
$\D$ be a block of $G_k $ with $|\D|>1$. Assume $2\in \D$. Suppose
$\D\subseteq \O_2$. Since $|\D|>1$, let $2\not=3\ell+2\in \O_2$. By
considering the permutation $(b_kb_k^{a_k})^{a_k}$, $\D$ contains at
least one element in $\O_1$, a contradiction. Thus, there exist
$x\in \O_1\cup\O_3$ and $x\in \D$. Since $b_kb_k^{a_k}=(e_k\ f_k)$
fixes $2$ and is cyclic on $\O_1\cup\O_3$, one has
$\O_1\cup\O_3\subseteq \D$, implying that $|\D|\geq 8k+8$. Since
$|\D|$ is a divisor of $12k+10$, one has $|\D|=12k+10$. Thus, $G_k$
is primitive on $\O$. By Proposition~\ref{Jordan}, $G_{10}=A_{10}$
and for $k\geq 1$, by Proposition~\ref{Marggraf}, $G_k$ is
$(4k+4)$-transitive because $G_k$ contains a cycle of length $8k+7$.
Since $4k+3>6$, $G_k=A_{12k+10}$. \hfill
\end{proof}

%
%

\section{Chains of consecutive $2$-covers}
\label{sec-height}

We start with some terminology. 
Let
$$%
X_n \stackrel{\p_n}{\to} X_{n-1}  \cdots X_2 \stackrel{\p_2}{\to}  X_1 \stackrel{\p_1}{\to} X_0
$$%
be a chain of  consecutive regular covering projections.
Suppose further that there is a chain of groups $G = G_0, G_1, G_2, \ldots, G_n$ such that 
$G_j \leq \Aut\, X_j$ is the lift of $G_{j-1} \leq \Aut\, X_{j-1}$ along $\p_j$ for each $j \in \{1,2,\ldots,n\}$. We then say that 
the above chain of covers is 
{\em $G$-admissible}. A  chain with this property  is denoted by  $\C(X,G)$. 
In particular, if all extensions $\CT(\p_j) \to G_j \to G_{j-1}$ are split, then the
chain is said to be {\em $G$-split-admissible}.
Further,  a {\em $G$-split-admissible}  chain $\C(X,G)$ is said to be
{\em split-sectional} and {\em split-transitive} provided all
extensions $\CT(\p_j) \to G_j \to G_{j-1}$
are,  respectively, split-sectional and split-transitive.

Let  $G=G_0$ be an arc-transitive group of automorphisms of a  symmetric  graph $X = X _0$. 
The {\em length} of the pair $(X,G)$ is the largest integer $n$ such that there exists a $G$-admissible 
chain of $n$ consecutive $2$-covers $X_n \to X_{n-1} \to  \ldots  \to X_1 \to X_0$. In particular, the 
{\em split-length} of the pair $(X,G)$ is the 
largest integer $n$ such that there exists a $G$-split-admissible 
chain of $n$ consecutive $2$-covers $X_n \to X_{n-1} \to  \ldots  \to X_1 \to X_0$.
Analogously,
the {\em sectional-length} and the {\em transitive-length}
of the pair $(X,G)$ is the 
largest integer $n$ such that there exists, respectively,
a sectional $G$-split-admissible chain and a transitive $G$-split-admissible chain  
of $n$ consecutive $2$-covers $X_n \to X_{n-1} \to  \ldots  \to X_1 \to X_0$.

Note that for a  {\em $G$-split-admissible} chain $\C(X,G)$
we have that $G_j\cong G \times \ZZ_2^j$ for
every $1 \le j\le n$.  Moreover,  since  $G_n$ is
arc-transitive on $X_n$  we have that  $G_n=\langle H_n, b\rangle$,
where $H_n$ is the stabilizer of a vertex  in $G_n$ and $b$ maps this
vertex to one of its neighbors. Since $G_n\cong G\times\mz_2^n$ has
at least $n$ generators, one has $n\le |H_n|$. Recall that $|H_n|$
is equal to the order of stabilizers in $G$. 
Hence, 

\begin{proposition}
\label{pro:finite}
The split-length of a symmetric  graph relative to an arc-transitive
group of automorphisms  must be finite.
\end{proposition}

\begin{lemma}
\label{lem:2step}

Let $G$ be an arc-transitive group of automorphisms of a symmetric  graph $X$, and let $\C(X,G)$ be 
a $G$-admissible chain 
$$%
X_2 \stackrel{\p_2}{\to} X_1 \stackrel{\p_1}{\to} X_0
$$%
of two consecutive $2$-covers where $\p_2$ is  split-transitive or split-sectional, respectively.
Then there exists  a $G$-admissible chain $\C'(X,G)$
$$%
X_2 \stackrel{\p_2'}{\to} X_1 '\stackrel{\p_1'}{\to} X_0
$$%
of two consecutive $2$-covers such that 
$\p_1'$ is  split-transitive or split-sectional, respectively.
\end{lemma}

\proof
Denote the respective lifted groups in the above sequence by 
$G = G_0, G_1, G_2$, and let  
$\CT(\p_1) = \la c_1 \ra$ and $\CT(\p_2) = \la c_2 \ra$.
Since $\p_2$ is split, we
have $G_2 = \bG_1 \times \la c_2\ra$, where $\bG_1 \cong G_1$. 
Since  $G_1$ contains $\la c_1 \ra  \cong \ZZ_2$ as a normal subgroup, 
$G_2$ contains $\la \bar{c}_1 \ra \times \la c_2 \ra \cong \ZZ_2 \times \ZZ_2$,
where $\la \bar{c}_1 \ra \leq \bG_1$ projects onto $\la c_1 \ra$ along $\p_2$. Of course, $\la \bar{c}_1 \ra$ is normal
in $G_2$  and hence 
acts semiregularly on $X_2$. Moreover, observe that the orbits of $\la \bar{c}_1\ra \times \la c_2\ra$ contain no edges.
 Therefore, 
quotienting by the action of  $\la \bar{c}_1 \ra$ gives  rise to a covering projection 
$\p_2' \colon X_2 \to X_1'$, where the quotient graph $X_1'$ is again simple. 
(Note: it might happen that $X_1' = X_1$.)
Let $\la \underline{c}\,_2 \ra$ be the projection of 
$\la c_2 \ra$ along $\p_2'$. Then the group 
$G_2$ projects onto 
$G_1' = \bG_1/\la \bar{c}_1 \ra \times \la \underline{c}\,_2\ra \cong G_0 \times \la\underline{c}\,_2\ra $.
Quotienting by $\la\underline{c}\,_2\ra$ gives rise to a covering projection 
$\p_1' \colon X_1' \to X_0$. 

If $\p_2$ is split-transitive, then  $\p_1'$ is  split with  a transitive complement since 
$\bG_1/\la \bar{c}_1 \ra$ is transitive on $X_1'$ which projects isomorphically onto $G_0$.
In fact,  $\p_1'$ is also split-transitive.
Suppose on the contrary that there exists a sectional complement $H$ to $\la\underline{c}\,_2\ra$ within
$G_1'$. Clearly, its lift $\bH$ along $\p_2'$ is  not transitive on $X_2$, and $G_2 = \bH \times \la c_2 \ra$.
Now $\bH$ projects onto $G_1$ along $\p_2$. Hence  $\p_2$ is split-mixed, a contradiction.

If $\p_2$ is split-sectional,  then  $\p_1'$ is  split since $\bG_1/\la \bar{c}_1 \ra$
is a complement to $\la\underline{c}\,_2\ra$ which projects onto $G_0$ along $\p_1'$.
We now show that $\bG_1/\la \bar{c}_1 \ra$  is indeed  sectional.  First note that by Proposition~\ref{niceCDC}
the projection $\p_2$ is the 
canonical double cover. Hence  $X_2$ is bipartite. Now the group $\bG_1$ preserves
the bipartition sets of $X_2$, which implies that $\bG_1/\la \bar{c}_1 \ra$ also preserves the 
bipartition sets of the graph $X_1'$.  Thus, $\p_1'$ is split with a sectional complement $\bG_1/\la \bar{c}_1 \ra$.
In fact,  $\p_1'$ is also split-sectional.
Suppose on the contrary that there exists a transitive complement $H$ to $\la\underline{c}\,_2\ra$ within
$G_1'$. Clearly, $G_2 = \bH \times \la c_2 \ra$ where $\bH$ is the lift of $H$ along $\p_2'$. Thus, 
$\bH$ is a sectional complement to $\la c_2 \ra$ in $G_2$ because $\p_2$ is assumed split-sectional.
Consequently, $\bH$ fixes the bipartition sets of $X_2$ setwise.
 Hence $\bH = \bG_1$ because   $G_2$ has 
a unique intransitive index $2$ subgroup. But then $H = \bG_1/\la \bar{c}_1 \ra$,  which is impossible 
 as $H$ is transitive on $X_1'$ and $\bG_1/\la \bar{c}_1 \ra$ is not. 
This completes the  proof.
\hfill\Qed

\bigskip
\noindent
{\bf Remark.}
As for the case when in Lemma~\ref{lem:2step} the projection $\p_2$ is  split-mixed,
from the  proof it is easy to see that 
there exist two chains 
\begin{center}
$X_2 \stackrel{\p_2'}{\to} X_1 '\stackrel{\p_1'}{\to} X_0$ \\
$X_2 \stackrel{\p_2''}{\to} X_1 ''\stackrel{\p_1''}{\to} X_0$
\end{center}
such that $\p_1'$ is split-transitive or split-mixed and  
$\p_1''$ is split-sectional or split-mixed. However, it may happen that neither 
$\p_1'$  nor  $\p_1''$ is split-mixed.
 See Example~\ref{ex:petersen}.

\bigskip
In what follows we shall restrict ourselves to cubic symmetric graphs.
Theorem~\ref{the:atmost2} below gives us some partial information about 
the structure of consecutive $2$-covers of such graphs.
As a consequence
we show  in Corollary~\ref{cor:split} that the split-length of a cubic symmetric graph relative to 
an  arc-transitive
group of automorphisms is at most $2$.

\begin{theorem}
\label{the:atmost2}

Let $G$ be a arc-transitive group of automorphisms of a symmetric cubic  graph $X$, and let $\C(X,G)$ be 
a $G$-admissible chain 
$$%
X_n \stackrel{\p_n}{\to} X_{n-1}  \cdots X_2 \stackrel{\p_2}{\to}  X_1 \stackrel{\p_1}{\to} X_0
$$%
of consecutive $2$-covers. Then at most two of $2$-covers in 
the chain are split. In particular, if  precisely two of them   are split, then one is split-transitive and the other is 
split-sectional or split-mixed, and this can happen only if $G$ is $s$-regular for $s = 2,3,5$.
\end{theorem}

\proof
Clearly, at most one of the $2$-covers in the chain $\C(X,G)$  can be split with a sectional complement
(that is, split-sectional or split-mixed). Namely,
suppose that $\p_i \colon X_i \to X_{i-1}$ is the first  such $2$-cover. 
Then by Proposition~\ref{niceCDC},  the projection
$\p_i$ is the  canonical double cover; hence all  $X_j$, $j \geq i$,   are bipartite.
Because graphs are assumed connected, all  $\p_j$, $j \geq i$, are split-transitive or non-split. 

Applying Lemma~\ref{lem:2step}, we may without  loss  of generality  assume  that all
split-transitive  covers come first followed by all others.
We now show that the initial subchain containing only split-transitive
 covers is of length at most $1$.
Suppose on the contrary that in the subchain
$$%
X_2 \stackrel{\p_2}{\to} X_1 \stackrel{\p_1}{\to} X_0
$$%
both $\p_1$ and  $\p_2$ are split-transitive.

Suppose that $G$ is $s$-regular.
First note that   $s \neq 1, 4$. Namely,  Theorem~\ref{transitive14} implies that each
split $2$-cover of $X_0$ relative to $G$ is  the  canonical double cover. Hence $X_1$ is bipartite.
But then $X_2$ is disconnected, a contradiction.
Also, note that in this case at most one split $2$-cover exists in $\C(X,G)$.

We may therefore restrict ourselves to the case  $s \in \{2,3,5\}$. 
In what follows we shall be using the following fact: a vertex-transitive index $2$ subgroup in an
$s$-regular group, where $s \geq 2$, must be $(s-1)$-regular.

Denote now the corresponding lift of $G$ along $\p_1$ by 
$G_1 =  \bG\times \la c_1 \ra $, where $\CT(\p_1) = \la c_1\ra$,  and 
the lift of $G_1$ along $\p_2$ by 
$G_2 =  \bG_1 \times  \la c_2 \ra =  \tilde{G} \times \la \bar{c}_1 \ra \times \la c_2 \ra$, where
$\CT(\p_2) = \la c_2\ra$, and $\bar{c}_1$ and $\tilde{G}$  project to  $c_1$ and $\bG$ along $\p_2$,
respectively.
 Then $G_1$ and $G_2$ are $s$-regular, and the transitive complements 
$\bG$ within $G_1$ and $\bG_1$ within $G_2$ are $(s-1)$-regular. Similarly,
the transitive complement $\bG_1 = \tilde{G} \times \la \bar{c}_1\ra$ to $\la c_2\ra$ within $G_2$, and the 
lift $\tilde{G} \times \la c_2 \ra$ of $\bG$ along $\p_2$, are $(s-1)$-regular. 
Suppose that $\tilde{G}$ is not transitive. Then the group $\bG$  
lifts along $\p_2$ to $\tilde{G} \times \la c_2 \ra$ with $\tilde{G}$ as a sectional complement.
Hence $\p_2$ is the  canonical double cover, by Proposition~\ref{niceCDC}. But then, again by Proposition~\ref{niceCDC},
the projection $\p_2$ cannot be split-transitive relative to $G_1$, a contradiction.
We conclude that $\tilde{G}$ is transitive on $X_2$.
This is an immediate contradiction if $s= 2$  because  $\tilde{G} \times \la c_2 \ra$ is then $1$-regular and cannot 
contain a transitive subgroup of index $2$. If $s \geq 3$, then $\tilde{G}$ is  $(s-2)$-regular. But that is again a contradiction:
if  $s=5$, we get a contradiction  because a $5$-regular group $G_2$ cannot contain a $3$-regular subgroup, by \cite{DM};
if $s = 3$ we get a contradiction because by \cite{CN},   a $3$-regular group $G_2$ cannot contain three $2$-regular subgroups
(namely, $\tilde{G} \times \la c_2 \ra$,
$\tilde{G} \times \la \bar{c}_1\ra$, and $\tilde{G} \times \la \bar{c}_1c_2\ra$). 

This shows that there is at most one split-transitive $2$-cover in the chain, and the proof is complete.
\hfill\Qed

\bigskip

Let $X$ be a cubic symmetric graph, and let $G$ be an arc-tarnsitive subgroup of $\Aut \, X$.
We say that $G$ is of type $(2^1)$ if it is  $2$-regular and if it
contains an involution flipping an edge.
Similarly, we say that $G$ is of type $(2^2)$ if it is  $2$-regular and if it
contains no involution flipping an edge.
If the full automorphism group $\Aut\,X$ is of type $(2^1)$ (respectively, of type $(2^2)$) and $\Aut\,X$ contains no other 
arc-transitive subgroups, then $X$ is called of type $(2^1)$ (respectively, of type $(2^2)$). 
Next, we say that $X$ is o type $(2^1,3)$ if $\Aut \, X$ is $3$-regular
and contains a $2$-regular subgroup of type $(2^1)$, but no other arc-transitive subgroups.
And finally, we say that $X$ is of types $(3)$ or  $(5)$, respectively,  if it is $3$-regular or  $5$-regular, respectively, 
and if, furthermore, $\Aut \, X$ contains no other arc-transitive subgroups.
With this terminology, we have the following corollary.

\begin{corollary}
\label{cor:split} 
Let $X$ be a connected cubic symmetric graph and let $G$ be an
$s$-regular subgroup of $\Aut \, X$. 
If $s \in \{1,4\}$ then the split-length  of $(X,G)$  is at most $1$.
If  $s \in \{2,3,5\}$, then  the split-length    is at most $2$, and moreover, if it is $2$ then
$X$ is  one of the following types: $(5)$, $(3)$, 
$(2^1)$,  or $(2^1,3)$ (with $G$ a $2$-regular subgroup of $\Aut \, X$). 
 \end{corollary}

\begin{proof}
Using Theorem~\ref{transitive14}
we have that $s \in \{2,3,5\}$ if the split-length of $(X,G)$ is $2$.
Furthermore, $G$ has no $1$-regular and no  $4$-regular subgroups. 
(If otherwise, Theorem~\ref{transitive14} implies that both covers in 
the chain  must be  canonical double covers, a contradiction.)

Suppose  that $G$ is of type $(2^2)$.
If there exists a split-transitive  $2$-cover $\p \colon \tilde{X} \to X$ relative to $G$, then
the complement $\bG \cong G$ of  $\CT(\p)$ within $\tilde{G} \cong G \times \ZZ_2$
is $1$-regular on  $\tilde{X}$. But this is a contradiction  since $\tilde{G}$ is also of type $(2^2)$, 
and by \cite[Theorem~3]{DM}  a group of type $(2^2)$ does not contain a $1$-regular subgroup.

Using  the list in \cite{CN} of all possible types of cubic symmetric graphs, this leaves us 
with $(5)$, $(3)$, $(2^1)$,  or  $(2^1,3)$  as the only possible types of $X$.
As for the case  when $X$ is of type $(2^1,3)$, suppose  that $G$ is $3$-regular.
If there exists a split-transitive $2$-cover $\p \colon \tilde{X} \to X$ relative to $G$, then
the complement $\bG \cong G$ of  $\CT(\p)$ within $\tilde{G} \cong G \times  \ZZ_2$
is $2$-regular on  $\tilde{X}$, and \cite[Theorem~4]{DM} implies that $G\cong H\times\ZZ_2$ where $H$
is a $2$-regular subgroup of $\Aut\,X$.  But then $X$ is a split-transitive 
$2$-cover of some smaller cubic symmetric graph,
a contradiction since by Theorem~\ref{the:atmost2} there is at most one split-transitive $2$-cover in the 
corresponding chain. This completes the proof of Corollary~\ref{cor:split}.
\end{proof}

%
%

\section{Concluding remarks}
\label{sec:remarks}
\indent

Let $X$ be a connected symmetric graph and let $G$ be an
$s$-arc-transitive subgroup of $\Aut \,X$. By
Proposition~\ref{pro:finite}, the split-length of $(X,G)$
is finite.  Furthermore, by Corollary~\ref{cor:split}, if $X$ is cubic,
then  the split-length of $(X,G)$ is at most $2$. 
This suggest the following  natural problem.

\begin{problem}
\label{prob:split}
Let $X$ be a symmetric graph with valency greater than $3$,  and let
$G$ be an arc-transitive subgroup of $\Aut \,X$. Find the split-length of  $(X,G)$.
\end{problem}

\begin{problem}
Let $X$ be a connected symmetric graph with valency greater than $2$, 
and let $G$ be an arc-transitive group of automorphisms of $X$. Can
the length of a chain of consecutive non-split $2$-covers relative
to $(X,G)$ be infinite? In particular, what can we say about the
case when $X$ is cubic?
\end{problem}

This brings us to (non-split) lengths of the pair $(X,G)$.
Let $C_{2n}$ be the cycle of length $2n$ with $n\geq 2$. Then, $\Aut\, C_{2n}$ is
arc-transitive and isomorphic to the dihedral group  $D_{4n}$. Clearly, $\Aut\, C_{2n}$
has a normal subgroup of order $2$ (the center of $\Aut\, C_{2n}$),
and so $C_{2n}$ is a $2$-cover of $C_n$ with $\Aut\, C_n$ lifting to
$\Aut\, C_{2n}$. If $n$ is even, this $2$-cover is non-split
relative to $\Aut\, C_n$,  for otherwise the center of $\Aut\,C_{2n}$ 
would contain  at least $4$ elements which is not the case.
Thus, for any even $n\geq 2$ the length of consecutive non-split
$2$-covers relative to $(C_n, \Aut\, C_n)$ is infinite.

As for cubic graphs, by
\cite{B88}, there is a unique cubic symmetric graph $\F40$ of order $40$ 
and  a unique such graph $\F80$ of order $80$.
Using the computer software package
MAGMA \cite{Mag}, one can  show that $\F80$ is a non-split $2$-cover of $\F40$
relative to $\Aut\, F40$, and that $\F40$ is a non-split $2$-cover of the
Desargues graph $\F20$B relative to $\Aut\, \F20$B. 
Since $\F20B$ is a split $2$-cover of the Petersen graph $\F10$, this
brings to length of a chain of consecutive $2$-covers to $3$.
We do not know, however, if such
a chain of length $4$ exists in the case of cubic symmetric graphs.
We would like to propose the following problem.

\begin{problem}
Let $X$ be a connected cubic symmetric graph, 
and let $G$ be an arc-transitive subgroup of $\Aut \, X$. 
Is there an upper bound on the length of the chain of consecutive $2$-covers relative
to $(X,G)$? 
\end{problem}

%
%

\begin{footnotesize}

\end{footnotesize}


\end{document}